\documentclass[11pt,a4paper,openany]{article}
\usepackage{mathrsfs}
\usepackage{bbm}
\usepackage{amsfonts}
\usepackage{amsmath,amsfonts,ams}
\usepackage{latexsym}
\usepackage{amssymb}
\DeclareMathOperator*{\esssup}{ess\; sup}

\numberwithin{equation}{section}
\newtheorem{proposition}{Proposition}[section]
\newtheorem{definition}{Definition}[section]
\newtheorem{lemma}{Lemma}[section]
\newtheorem{theorem}{Theorem}[section]

\newtheorem{Rem}{Remark}[section]

\begin{document}
\title{On the Real Analyticity of the Scattering Operator  for  the Hartree  Equation  }

\author{Changxing Miao$^1$, Haigen Wu $^{2, 3}$ and  Junyong Zhang $^2$\\
        \\
{\small $^{1}$Institute of Applied Physics and Computational Mathematics}\\
        {\small P. O. Box 8009, Beijing, China, 100088 } \\
        {\small (miao\_changxing@iapcm.ac.cn )}        \\
       {\small $^{2}$The Graduate School of China Academy of Engineering Physics  }\\
        {\small P. O. Box 2101,\ Beijing,\ China,\ 100088 ;} \\
{\small $^{3}$School of Mathematics and Information Science}\\
{\small Henan Polytechnic University, Jiaozuo, China, 454000} \\
        {\small ( wuhaigen@hpu.edu.cn, zhangjunyong111@sohu.com )}}
\date{}
\maketitle

\maketitle
\begin{abstract}
In this paper,  we study the real analyticity of the scattering
operator for the
 Hartree  equation $ i\partial_tu=-\Delta u+u(V*|u|^2)$.
To this end,  we  exploit interior and exterior cut-off in time
and space,  and combining with the compactness argument to
overcome difficulties which arise from absence  of good properties
for the nonlinear Klein-Gordon equation, such as the finite speed
of propagation and ideal time decay estimate. Additionally, the
method in this paper allows us to simplify the proof of
analyticity of the scattering operator for the
 nonlinear Klein-Gordon equation with cubic nonlinearity in Kumlin\cite{Ku}.
\end{abstract}



\section{Introduction}
This paper is devoted to the proof of the real analyticity of
scattering operator for the Hartree equation
\begin{equation}\label{1.1}
\begin{cases}
i\partial_tu=-\Delta u+u(V*|u|^2),~~~\quad  (t, x)\;\in\; \R\times\R^3,\\
u(0)=u_0(x)\in H^1(\R^3).\end{cases}
\end{equation}
Here $u(t,x)$ is a complex valued function defined in $\R^{1+3}$, $V(x)$, called potential,
is a real valued radial function
defined in $\R^3$, and $*$ denotes the convolution in $\R^3$.
Under suitable assumption on $V$, Ginibre-Velo\cite{GV} proved the scattering theory of the equation \eqref{1.1} in the
energy space $H^1$. Attempting  to study the (complex) analyticity
of the scattering operator is in vain  because $\bar{u}$ is not
analytic even if $u$ is. However, following the W. Strauss  suggestion (private communication), we
can study the real analyticity which is still a very interesting issue. \vskip0.1cm

Let $~u=\varphi (t,x)+i\psi(t,x),~u_0=\varphi _0(x)+i\psi_0(x)$,
and$~\varphi (t,x),~\psi(t,x),~\varphi _0(x),~\psi_0(x)~$ are real
valued functions defined in $\R\times\R^3$ or $\R^3$. Then the integral
form of equation\eqref{1.1}
\begin{align}
u(t)=e^{it\Delta}u_0-i\int_0^t
e^{i(t-s)\Delta}((V*|u|^2)u)(s)ds\label{1.2}
\end{align}can be rewritten as
\begin{align}
\binom{\varphi }{\psi}=&\begin{pmatrix}\cos t\Delta &-\sin t\Delta \\
\sin t\Delta & \cos
t\Delta\end{pmatrix}\binom{\varphi _0}{\psi_0}\nonumber\\
&+\int_0^t\begin{pmatrix} \sin (t-s)\Delta & \cos (t-s)\Delta\\-\cos
(t-s)\Delta &\sin (t-s)\Delta\end{pmatrix}\binom{\varphi
}{\psi}\big(V*(\varphi ^2+\psi^2)\big)(s)ds.\label{1.3}\end{align}
Setting
$$U(t)=\binom{\varphi (t)}{\psi(t)}~~\text{and}~~ U_0=\binom{\varphi_0}{\psi_0},$$
then \eqref{1.3} can be transformed into
\begin{align}
\mathcal {N}(t)U_0:=U(t)=G(t)U_0-\int_0^t
\Delta^{-1}G'(t-s)(V*|U|^2)U(s)ds,\label{1.4}
\end{align}
where $$G(t)=\begin{pmatrix}\cos t\Delta &-\sin t\Delta \\
\sin t\Delta & \cos t\Delta\end{pmatrix}$$ is  a unitary  group
associated with the equation\eqref{1.1}.

\medskip
First, we recall the  decay estimate and Strichartz estimates in the
context of Schr\"odinger equation( see
\cite{Ca},\cite{KT},\cite{Mi}).
\begin{definition} A pair $(q,r)$ is admissible, denoted by $(q, r)\in \Lambda$,  if $r\in [2,6]$ and
$q$ satisfies
\begin{align}\frac{2}{q}=\delta(r):=3\Big(\frac{1}{2}-\frac{1}{r}\Big). \end{align}
\end{definition}

\begin{lemma}\label{lem1.1}Let $S(t)=e^{it\Delta}$, then

(1)the $L^{r'}-L^r$ decay estimate
\begin{align} \|S(t)\varphi \|_r\leq C|t|^{-\delta(r)} \|\varphi(x) \|_{r'},\label{1.5}
\end{align}holds for $2\leq r\leq\infty $ ;

(2) the Strichartz estimates
\begin{align} \|S(t)u \|_{L^q(\R,L^r(\R^3))}&\leq C \|u \|_2,\label{1.6}\\
\Big\|\int_0^tS(t-s)f(s)ds\Big\|_{L^{q_1}(I,L^{r_1}(\R^3))}&\leq
C \|f \|_{L^{q'_2}(I,L^{r'_2}(\R^3))} \label{1.7}
\end{align}hold  true for any interval $I\subset \R$, and  for any admissible pairs  $(q,r),~(q_j,r_j)\in \Lambda,~j=1,2$.
\end{lemma}
\begin{Rem}\label{rem1.1} Lemma \ref{lem1.1} still holds for the propagators
$G(t)$ and $\Delta^{-1}G'(t)$ by Euler formulae $$\cos
t\Delta=\frac{e^{it\Delta}+e^{-it\Delta}}{2},~\sin
t\Delta=\frac{e^{it\Delta}-e^{-it\Delta}}{2}.$$
\end{Rem}

Let  $B$ be  a Banach space, and
 $$\binom{u}{v}\in B \quad\Longleftrightarrow\quad  \Big\|\binom{u}{v}\Big\|_B= \|u \|_B+ \|v \|_B<\infty.$$
Throughout this
paper, the symbol  $C$ denotes a constant which may be
different from line to line, and $C(*)$ denotes the constant which only
depends on the parameter $*$.

\smallskip
 Define the wave operator $W_{\pm}:U_{\pm}\mapsto U_0$
as follows: for any $U_{\pm}\in H^1$, there exists $U_0\in H^1$ such
that
\begin{align} \label{1.8}\big\|G(t)U_{\pm}-\mathcal
{N}(t)U_0\big\|_{H^1}\longrightarrow 0,\quad\text{as~}
t\longrightarrow\pm\infty.
\end{align}
When  the wave operator $W_{\pm}$ are invertible operators,  we can define scattering  operator as
$S=W_+^{-1}\circ W_-:U_-\mapsto U_+$.

Set
$$X=C(\R,H^1(\R^3))\cap
\bigcap_{(q,r)\in\Lambda}L^q(\R,H^1_r(\R^3)).$$
Ginibre-Velo established a complete scattering theory in energy space  provided that
the potential $V$ satisfies the
following assumption:

\begin{enumerate}
\item[(H1)] $V$ is a real function and $V\in
L^{p_1}+L^{p_2}$ for some $p_1, p_2$ satisfying
$$1< p_2 \leq p_1 <\frac32.$$
\end{enumerate}

\begin{enumerate}
\item[(H2)] $V$ is radial and nonincreasing, namely $V(x)=v(r)$
where $v$ is nonincreasing in $\mathbb{R}^+$. Furthmore, for some
$\alpha \geq 2, $ $v$ satisfies the following condition:

\vskip0.3cm
 ($A_{\alpha}$):\qquad There exists $a>0$ and $A_{\alpha}>0$ such that
$$v(r_1)-v(r_2) \geq \frac{A_{\alpha}}{\alpha} (r^{\alpha}_2 -r^{\alpha}_1)\ \text{for}\ 0 < r_1<r_2 \leq a.$$
\end{enumerate}
In particular,  the wave operator $W_{\pm}$ and the scattering operator $S$ are
bounded and continuous from $H^1$ to $H^1$.
\medskip

 Our main result is
\begin{theorem}\label{th1} Let $V(x)$ satisfy the assumption (H1) and (H2).
Then the operators $W_{\pm}$ and $S$ are analytic from $H^1$ to $H^1$.
\end{theorem}
The proof of Theorem \ref{th1} depends on the following theorem:
\begin{theorem}\label{th2}
Let $U_0\in H^1$ and $U(t)$ be the unique solution of \eqref{1.4}
in $X$. Then the map $U:U_0 \mapsto U(U_0)$ is analytic from $H^1$
to $X$.
\end{theorem}

\medskip
For the nonlinear  Klein-Gordon equation with cubic nonlinearity, using the contraction mapping
principle, Baez-Zhou\cite{BZ1} proved  the analyticity of
scattering operator on a neighborhood of the space of
finite-energy Cauchy data, $H^1\oplus L^2(\R^3)$. Kumlin\cite{Ku}
generalized the result to entire energy space by means of the
Fredholm alternative theorem. The proof in \cite{Ku} depends on
the following two good properties of the linear Klein-Gordon equation:
\vskip0.2cm

 (1) $L^p-L^{p'}$ estimates stated in the following proposition.
 \begin{proposition}[\cite{Br}\cite{Ku}]  Let $K(t)=\frac{\sin t(-\Delta+m^2)^{\frac{1}{2}}}{(-\Delta+m^2)^\frac{1}{2}}$,
 $1<p\leq 2\leq p',~\frac{1}{p}+\frac{1}{p'}=1$ and $\sigma:=\frac{1}{2}-\frac{1}{p'}$ and
 $0\leq \theta\leq 1$. Then if $ (n+1+\theta)\sigma\leq 1+s-s',$
 \begin{align}\label{k}\|K(t)g\|_{W^{s',p'}}\leq k(t)\|g\|_{W^{s,p}},~~t\geq 0,\end{align}
 where $$k(t)=\begin{cases}&t^{-  (n-1-\theta)\sigma},~0<t<1,\\
 &t^{-(n-1+\theta)\sigma},~1\leq t.\end{cases}.$$
 \end{proposition}
For suitable  $p,~p'$ and $\theta$ , $k(t)\in L^1(\R)$. In particular,
$$\|K(t)g\|_{H^1}\leq C\|g\|_{L^2}.$$
These estimates are crucial in the proofs of Step 1 and Step 3
in \cite{Ku}.
\vskip0.2cm

 (2) The finite speed of propagation.
\vskip0.1cm

The finite speed of propagation of the solution of linear wave equation means that  for $t\in
[-T,T],~T<\infty$,
\begin{align}\label{speed} \Big\|\int_0^tK(t-s)u(s)ds\Big\|_{L^6(|x|>R)}\thicksim
\Big\|\int_0^tK(t-s)\eta_Ru(s)ds\Big\|_{L^6(\R^3)}. \end{align}
Namely, the cut-off function $\eta_R(x)$ defined below in
\eqref{cut} commutes with the group $K(t)$ in some sense, which
plays an important role in the proof of analyticity, see  Step 2
in \cite{Ku}.

\vskip0.1cm
 The arguments in this paper still take advantage of  the Fredholm
alternative theorem  together with the analytic version of implicit
 function theorem (cf \cite{BZ1,BZ2, Ku}). However we have to overcome
 some difficulties arise from loss of the
good properties (1) and (2) for the Schr\"{o}dinger equation. Our
major innovations are as follows : Comparing with $k(t)$ in
\eqref{k}, the kernel $|t|^{-\delta(r)}$ in \eqref{1.5} is not in
any $L^p(\R),~1\leq p\leq\infty$, and  the
Hardy-Litttlewood-Sobolev inequality can not supply any decay yet.
A new approach to deal with the singular kernel is the double
localization in time
\begin{align}|t-s|^{-\delta(r)}\chi_{\{|s|\leq
T/2\}}\chi_{\{|t|>T\}},~~ |t-s|^{-\delta(r)}\chi_{\{|s|>
T/2\}}\chi_{\{|t|>T\}},
\end{align}
this together with other techniques
 helps us to get time decay,  where
 $\chi_{A}$ denotes the characteristic function on the interval $A$.
On the other hand, replacing \eqref{speed} by introducing  double (interior and
exterior) cut-off on space
\begin{align}
&\Big\|\int_0^t\Delta^{-1}G'(t-s)\big((V*|U|^2)\xi_R\Psi_j\big)(s)ds\Big\|_{L^4(|t|\leq
T,H^1_3(|x|>M))}~,\\
&\Big\|\int_0^t\Delta^{-1}G'(t-s)\big((V*|U|^2)\eta_R\Psi_j\big)(s)ds\Big\|_{L^4(|t|\leq
T,H^1_3(\R^3))}~,\end{align}
we obtain decay estimates and overcome the above difficulties by means of
 compactness principle,  where
 \begin{align}\label{cut}&\xi_R(x)  \in
 C_0^\infty(\R^3), \nonumber\\
&\xi_R(x)=1\text{~if~} |x|\leq R,~\xi_R(x)=0 \text{~if~}|x|>2R\quad \text{~and~}\nonumber\\
& \eta_R(x)=1-\xi_R(x). \end{align}
Based on these decay estimates and the arguments in  Kumlin\cite{Ku},
 we can prove  Theorem \ref{th1} by  the approximate theorem of analytic operator sequence (cf.\cite{HP}).
 However, it is worth mentioning that we  make  use of compactness arguments  and
 the definition of Frech\'{e}t
derivative to avoid repeating the argument of global time space
integrability,  and give a more concise proof of Theorem \ref{th1}.

  The paper is organized as follows: In Section 2 we give the proof of
  Theorem \ref{th2}, Section 3 is devoted to the  proof of key lemma which consists of the main part
   of this paper.  At last,  we  supply a brief derivation  of Theorem \ref{th1} in Section 4.

\section{Proof of Theorem \ref{th2}}
 To prove Theorem \ref{th2}, we need the following analytic version of implicit
 function theorem.

\begin{lemma}$\cite{BZ2}$\label{lem2.1}  Suppose that
$X,Y,Z$ are Bananch spaces and $Q$ is an open neighborhood of the
point $(x,y)\in X\times Y$. Suppose that $f: Q\mapsto Z$ is
analytic, $f(x,y)=0$ and $D_2f(x,y):X\mapsto Z$ has a left inverse,
where $D_2$ indicates the  Frech\'{e}t derivative with respect to the second
variable. Then for some open set $P$ containing $x$,  there exists a
unique analytic function $g:P\mapsto Y$ such that $g(x)=y$ and
$f(x',g(x'))=0$ for all $x'\in P$.
\end{lemma}

\smallskip
For $U_0\in H^1$ and $\Psi\in X$, we consider the following mapping
\begin{align}
\label{2.1}
R(U_0,\Psi)=G(t)U_0-\int_0^t\Delta^{-1}G'(t-s)\big((V*|\Psi|^2)\Psi\big)(s)
ds-\Psi(t).
\end{align}
Since it is linear in
$U_0$ and multilinear in $\Psi$, we  know that  $R:H^1\times X\mapsto X$
 is analytic by  the nonlinear estimate in \cite{GV}
$$\bigg\|\int_0^t\Delta^{-1}G'(t-s)\big((V*|\Psi|^2)\Psi\big)(s) ds\bigg\|_X<\infty, \qquad \forall \; \Psi\in
X.$$
  On the other hand, $R(U_0,U)=0$ by
$\eqref{1.4}$. Hence it suffices to prove the invertibility of
\begin{align}D_2R(U_0,U):\; X\mapsto X \nonumber
\end{align}
for each $U_0\in H^1$. By the open mapping theorem, we
only need to prove that $D_2R(U_0,U)$ is injective and surjective.

For $U_0\in H^1,~\Psi\in X$, one  has
\begin{align}
\label{2.2} D_2R(U_0,U)(\Psi)(t)=&-2\int_0^t\Delta^{-1}G'(t-s)\big(V*(\Psi
U)\big)\Psi(s) ds\nonumber\\&-\int_0^t\Delta^{-1}G'(t-s)(V*|U|^2)\Psi(s)
ds-\Psi(t).
\end{align}

(1)\; The injectivity of $D_2R(U_0,U)$. \vskip0.3cm

For simplicity, we always assume that $V(x)\in L^p$.
 Let
$D_2R(U_0,U)\Psi=0$, then
\begin{align}
\label{2.3} \|\Psi \|_r&\leq
C\int_0^t|t-s|^{-\delta(r)}\Big( \big\|\big(V*(U\Psi)\big)U(s) \big\|_{r'}
+\big \|(V*|U|^2)\Psi(s) \big\|_{r'}\Big)ds\nonumber\\
&\leq 2C\int_0^t|t-s|^{-\delta(r)} \|V\ \|_p \|U \|^2_l \|\Psi
\|_rds,
\end{align} where $2=\frac{1}{p}+\frac{2}{l}+\frac{2}{r}$.

For every $p\in(1,\frac{3}{2})$, we can take $l=r\in(3,4)$ such
that
\begin{align}
\label{2.4} \|\Psi \|_r\leq C \|V \|_p \|U \|^2_{L^\infty
H^1}\int_0^t|t-s|^{-\delta(r)} \|\Psi \|_rds.
\end{align}
For each $t\in (0,T)$, one easily verifies that  by \eqref{2.4}
\begin{align}
 \|\Psi(t) \|_r\leq Ct^{1-\delta{(r)}}\esssup_{s\in (0,t)} \|\Psi(s) \|_r.\nonumber
\end{align}  We choose $T$  small enough such that \begin{align}
\esssup_{t\in (0,T)} \|\Psi(t) \|_r \leq \frac{1}{2}\esssup_{s\in
(0,T)} \|\Psi(s) \|_r.\nonumber
\end{align}
This implies that  $\Psi \equiv 0,$ a.e. $t\in [0,T]$ for some $T>0$.
Repeating this process on $(nT,nT+T), ~n\in\Z$, we have $\Psi
\equiv 0,$ a.e. $(x,t)\in \R^{3+1}$.
\vskip0.3cm

(2)\;  The surjectivity of $D_2R(U_0,U)$. \vskip0.3cm

Setting\begin{align} \label{2.5}\mathcal
{T}_{U_0}\Psi(t)=&-\int_0^t\Delta^{-1}G'(t-s)\big((V*|U|^2)\Psi\big)(s)ds\nonumber\\
&-2\int_0^t\Delta^{-1}G'(t-s)\Big(\big(V*(U\Psi)\big)U\Big)(s)ds,
\end{align}
we have $D_2R(U_0,U)=\mathcal {T}_{U_0}-I$.   By the Fredholm
alternative theorem, our first choosing is  to show that $\mathcal {T}_{U_0}$
is compact operator from $X$ to $X$ since $\mathcal {T}_{U_0}-I$ is
injective. However, $\mathcal {T}_{U_0}$ may be not compact.
To our goal, it suffices to show that $\mathcal {T}_{U_0}^2$ is compact. In fact, since
\begin{align}\label{2.6}
\mathcal {T}_{U_0}^2-I=(\mathcal {T}_{U_0}-I)(\mathcal {T}_{U_0}+I)
\end{align}
and $\mathcal {T}_{U_0}+I$ is also injective, the
Fredholm theorem still works. Therefore \eqref{2.6} implies that the surjectivity
of $\mathcal {T}_{U_0}-I$.

Concerning the trilinear form
\begin{align}
\label{2.7} \mathcal
{B}(\Psi_1,\Psi_2,\Psi_3):=-\int_0^t\Delta^{-1}G'(t-s)\big(V*(\Psi_1\Psi_2)\big)\Psi_3(s)ds,
\end{align}
we have the following nonlinear estimate.

\begin{lemma}\label{lem2.2}
For $\Psi_j\in X,~j=1,2,3$ , one has
\begin{align}
\label{2.8} \big\|\mathcal {B}(\Psi_1,\Psi_2,\Psi_3)\big\|_X\lesssim
\prod_{j=1,2,3}
 \|\Psi_j \|_{L^4H_3^1}.\end{align}\end{lemma}

 {\bf Proof } \; Using the Strichartz estimates together with the H\"{o}lder inequality,  we obtain
\begin{align}\label{2.9}
&\Big\|\int_0^t\Delta^{-1}G'(t-s)\big(V*(\Psi_1\Psi_2)\big)\Psi_3(s)ds\Big\|_X\nonumber\\\lesssim
&
\big\|\big(V*(\Psi_1\Psi_2)\big)\Psi_3\big\|_{L^{4/3}H^1_{3/2}}\nonumber\\
\leq&\sum_{\{i,j,k\}=\{1,2,3\}}
 \|V \|_p \|\Psi_i \|_{L^4H^1_3} \|\Psi_j \|_{L^4L^{\tilde{p}}} \|\Psi_k \|_{L^4L^{\tilde{p}}}\nonumber\\
\lesssim& \prod_{j=1,2,3}
 \|\Psi_j \|_{L^4H_3^1},
\end{align}where $\tilde{p}=\frac{6p}{4p-3},$ and we have used the embedding relation
$H^1_3\hookrightarrow L^{\tilde{p}}$.
\medskip

As a direct consequence of Lemma \ref{lem2.2}, we have
 $$\big\|\mathcal {T}_{U_0}\Psi(t)\big\|_X\lesssim  \|\Psi \|_{L^4H_3^1}.$$
 This implies that $\mathcal {T}_{U_0}: L^4H_3^1\mapsto X$ is bounded.
 Since the composition of compact operator and bounded operator is still compact,
it is enough to verify following key lemma.
\begin{lemma}\label{main} \;  Let $U_0\in H^1$, Then
\begin{align}\mathcal{T}_{U_0}:~~X\mapsto ~L^4H^1_3 ~~~\text{
is compact}.\nonumber
\end{align}\end{lemma}

\section{Proof of Lemma \ref{main}}

Let $\{\Psi_j\}_{j=0}^\infty$  be uniformly bounded in $X$, i.e. $
\|\Psi_j\|_X\leq C$ for constant $C>0$, we shall show that
$\{\mathcal {T}_{U_0}\Psi_j\}_{j=0}^\infty$ has a Cauchy
subsequence in $L^4H^1_3$. Our main tool is the Arzela-Ascoli
theorem, so it is necessary to localize both the time and the
space.   The proof can be divided into five steps.

Let$~\xi_R,~\eta_R ~\text{be defined as those in \eqref{cut}}$.
\begin{align}\text{Step 1.}
~~\lim_{T\rightarrow\infty}\sup_{j\in\N}\big\|\mathcal
{T}_{U_0}\Psi_j\big\|_{L^4(|t|>T,H^1_3(\R^3))}=0;~~~~~~~~~~~~~~~~~~~~~~~~~~~~~~~~~~~~~~~~~~~~
\nonumber\end{align}
\begin{align} \text{Step
2.}~~\lim_{R\rightarrow\infty}\sup_{j\in\N}\big\|\mathcal
{T}_{U_0}(\eta_R\Psi_j)\big\|_{L^4(|t|\leq T,H^1_3(\R^3))}=0
\text{~~for all~~} T>0,~~~~~~~~~~~~~~~~~~ ~~~  \nonumber\end{align}
\begin{align}
&\text{Step 3.}~~\lim_{M\rightarrow\infty}\sup_{j\in\N}\big\|\mathcal
{T}_{U_0}(\xi_R\Psi_j)\big\|_{L^4(|t|\leq T,H^1_3(|x|>M))}=0
\text{~for~all~} T>0,~R>0;~~~~~~~  \nonumber\end{align}

 Step 4.~~$\{\mathcal
{T}_{U_0}(\xi_R\Psi_j)\}_{j=0}^\infty$ has a Cauchy subsequence in
${L^4(|t|\leq T,H^1_3(|x|\leq
 M))}$ for all $T>0, ~R>0,~M>0;$

\smallskip
Step 5. A Cantor diagonalized  process.

\bigskip

For the sake of convenience,  we first give some useful estimates.
\begin{lemma}\label{lem3.1}
Let $ \|U \|_X\leq C$, $ \|\Psi_j \|_X\leq C$ and $V\in L^p(\R^3)$.
For any $p\in(1,\frac{3}{2})$, then
$$\big\|(V*|U|^2)\Psi_j\big\|_2\leq C$$ and
$$\big\|(V*|U|^2)\Psi_j\big\|_{\dot{H}^\delta}\leq C $$ hold for
sufficient small $\delta>0$. \end{lemma}

 {\bf Proof~} For any $p\in(1,\frac{3}{2})$,  by Sobolev embedding theorem it is derived that
 $$\big\|(V*|U|^2)\Psi_j\big\|_2\leq  \|V \|_p \|U\|_r^2 \|\Psi_j \|_r
  \leq \|V \|_p \|U \|^2_{L^\infty H^1} \|\Psi_j \|_{L^\infty H^1},$$
 where $1+\frac{1}{2}=\frac{1}{p}+\frac{3}{r},~ r\in(\frac{18}{5},6).$\\

For $p,r$ as above, taking $\delta>0$ small enough and  using
fractional Leibniz formula,   we have
\begin{align}&\big\|(V*|U|^2)\Psi_j\big\|_{\dot{H}^\delta}\nonumber\\\leq&
C( \|V \|_p \|U \|_r^2 \|D^\delta\Psi_j \|_r+ \|V \|_p \|U \|_r \|D^\delta
U \|_r \|\Psi_j \|_r)\nonumber\\
 \leq& \|V \|_p \|U \|^2_{L^\infty H^1} \|\Psi_j \|_{L^\infty H^1}.\nonumber\end{align}
This shows that the sequence $\{(V*|U|^2)\Psi_j\}_{j=0}^\infty $ is
uniformly bounded in  $H^\delta(\R^3)$.

\smallskip
\begin{lemma} \label{lem3.2} Let $T<\infty,~t\in [-T,T]$, then
\begin{align}\label{3.1}\Big\|\int_0^t\Delta^{-1}G'(t-s)U(s)ds\Big\|_{H^{2-\epsilon}(\R^3)}\leq
C \|U\|_{L^2([-T,T],L^2(\R^3))}
\end{align}holds for all $\epsilon>0$.\end{lemma}

{\bf Proof~}\;(cf.\cite{Br2}\cite{Ku})  Let $k\in\R$. $F^k$ denotes the
operator on $L^2[-T,T]$ defined by
$\widehat{F^kh}(n)=(in)^k\hat{h}(n),~n\in\N-\{0\},$ and
$\widehat{F^kh}(0)=\hat{h}(0)$, where $ h\in L^2[-T,T]$ and $\hat{}$ denotes the Fourier
transform. Making use of the discrete Plancherel identity and the
transformation between time and space regularity, it follows that by taking $2k=2-\epsilon$ with $k<1$
\begin{align}&\Big\|\int_0^t\Delta^{-1}G'(t-s)U(s)ds\Big\|_{H^{2-\epsilon}}\nonumber\\
\leq&\Big\|\int_{-T}^{T}\big(F^k\chi_{[0,t]}\big)(s)\cdot\big(F^{-k}\Delta^{-1}G'(t-\cdot)U(\cdot)\big)(s)ds\Big\|_{H^{2-\epsilon}}\nonumber\\
\leq& \|F^k\chi_{[0,t]} \|_{L^2([-T,T])}\cdot \|F^{-k}\Delta^{-1}G'(t-\cdot)U(\cdot) \|_{L^2([-T,T],H^{2-\epsilon})}\nonumber\\
\leq&C \|U \|_{L^2([-T,T],L^2)},\nonumber\end{align}
where we have used the following estimate
$$ \|F^k\chi_{[0,t]} \|_{L^2([-T,T])}\leq \Big(\sum_{n\in
\Z-\{0\}}(n^k\cdot\frac{1}{n})^2+1\Big)^\frac{1}{2}<\infty.$$

\medskip
Now, we are in position to prove Lemma \ref{main}.  Note that the
multi-linear estimates of two terms of $\mathcal{T}_{U_0}$ are
similar, we only need to estimate the first term.

 \medskip
{\bf Step 1~}   We make use of  an
 interior time cut-off technique to deal with the convolution kernel. It is easy to show that

\begin{align}&\Big\|\int_0^t\Delta^{-1}G'(t-s)(V*|U|^2)\Psi_j(s)ds\Big\|_{L^4(|t|>T,H^1_3(\R^3))}\nonumber\\
\leq&\Big\|\int_0^t|t-s|^{-\frac{1}{2}} \|\big((V*|U|^2)\Psi_j\big)(s) \|_{H^1_{3/2}}ds\Big\|_{L^4(|t|>T)}\nonumber\\
\leq&\Big\|\int_0^t|t-s|^{-\frac{1}{2}}\chi_{\{|s|>\frac{T}{2}\}}(s)
\|\big((V*|U|^2)\Psi_j\big)(s)
\|_{H^1_{3/2}}ds\Big\|_{L^4(|t|>T)}\nonumber\\
&+\Big\|\int_0^t|t-s|^{-\frac{1}{2}}\chi_{\{|s|\leq\frac{T}{2}\}}(s)
\|\big((V*|U|^2)\Psi_j\big)(s)
\|_{H^1_{3/2}}ds\Big\|_{L^4(|t|>T)}\nonumber\\=&:I_1+I_2
\end{align}

On the one hand,
\begin{align}I_1=&\Big\|\int_0^t|t-s|^{-\frac{1}{2}}
\|\big((V*|U|^2)\Psi_j\big)(s)\chi_{\{|s|>\frac{T}{2}\}}(s) \|_{H^1_{3/2}}ds\Big\|_{L^4(|t|>T)}\nonumber\\
\lesssim& \|(V*|U|^2)\chi_{\{|s|>\frac{T}{2}\}}\Psi_j
\|_{L^{4/3}(\R,H^1_{3/2})}\nonumber\\
\leq& \|(V*|U|^2)\chi_{\{|s|>\frac{T}{2}\}} \|_{L^2(\R,L^{3})} \|\Psi_j \|_{L^4(\R,H^1_3)}\nonumber\\
&+ \|(V*|U|^2)\chi_{\{|s|>\frac{T}{2}\}} \|_{L^2(\R,H^1_{2p})} \|\Psi_j \|_{L^4(\R,L^{\tilde{p}})}\nonumber\\
\leq&\|\Psi_j \|_X\Big( \|(V*|U|^2)\chi_{\{|s|>\frac{T}{2}\}}
\|_{L^2(\R,L^{3})}+\|(V*|U|^2)\chi_{\{|s|>\frac{T}{2}\}}
\|_{L^2(\R,H^1_{2p})}\Big).\nonumber
\end{align}
Since
\begin{align}&\|(V*|U|^2)\chi_{\{|s|>\frac{T}{2}\}}\|_{L^2(\R,L^{3})}+\|(V*|U|^2)\chi_{\{|s|>\frac{T}{2}\}} \|_{L^2(\R,H^1_{2p})}\nonumber\\
\leq&
\|V\|_p\|U\|^2_{L^4(\R,L^{\tilde{p}})}+\|V\|_p\|U\|_{L^4(\R,L^{\tilde{p}})}\|U\|_{L^4(\R,H^1_3)}\nonumber\\
\leq& 2\|V \|_p \|U
\|^2_{L^4(\R,H^1_{3})}<\infty,\nonumber\end{align} we get
$$\lim_{T\rightarrow\infty} I_1=0$$
by
\begin{align} &~~~~\|(V*|U|^2)\chi_{\{|s|>\frac{T}{2}\}} \|_{L^2(\R,L^{3})}+ \|(V*|U|^2)\chi_{\{|s|>\frac{T}{2}\}}
\|_{L^2(\R,H^1_{2p})}\nonumber\\&= \|(V*|U|^2)
\|_{L^2(|t|>\frac{T}2,L^{3})}+ \|(V*|U|^2)
\|_{L^2(|t|>\frac{T}2,H^1_{2p} )}\nonumber\\
&\longrightarrow  ~0,~~\text{as~} T\longrightarrow~
+\infty,\end{align}
where we have used the property of absolute continuity.

On the other hand, similar arguments as deriving the estimate of $I_1$ can be used to get that
\begin{align}I_2&\leq
\Big(\int_{|t|>T}\Big|\int_\R|t-s|^{-\frac{1}{2}}
\chi_{\{|s|\leq\frac{T}{2}\}}(s)\|\big((V*|U|^2\Psi_j\big)(s)\|_{H^1_{3/2}}ds\Big|^4dt\Big)^\frac{1}{4}\nonumber\\
& \lesssim \Big(\int_T^\infty |t-\frac{T}{2}|^{-\frac{1}{2}\cdot
4}dt\Big)^\frac{1}{4}
\Big(\Big|\int_\R\chi_{\{|s|\leq\frac{T}{2}\}}(s)\|\big((V*|U|^2\Psi_j\big)(s)\|_{H^1_{3/2}}ds\Big|^4\Big)^\frac{1}{4}\nonumber\\
&\lesssim T^{-1/4}
\Big(\int_\R\chi_{\{|s|\leq\frac{T}{2}\}}(s)ds\Big)^{\frac{1}{4+\varepsilon}}
\Big(\int_\R
\|\big((V*|U|^2\Psi_j\big)(s)\|^\frac{4+\varepsilon}{3+\varepsilon}_{H^1_{3/2}}ds
\Big)^\frac{ 3+\varepsilon }{4+\varepsilon} \nonumber\\
&\lesssim T^{-{\frac{\varepsilon}{4(4+\varepsilon)}}}
\|\big((V*|U|^2)\Psi_j\big)(s)\|
_{L^\frac{4+\varepsilon}{3+\varepsilon}H^1_{3/2}}
  \nonumber\\
&\leq
T^{-{\frac{\varepsilon}{4(4+\varepsilon)}}}\Big(\|V\|_p\|\Psi_j\|_{L^4H^1_3}\|U\|^2_{L^qL^{r}}
+\|V\|_p\|U\|_{L^4H^1_3}\|U\|_{L^qL^{r}}\|\Psi_j\|_{L^qL^{r}}\Big),
\end{align}
here
$q=4-\frac{4\varepsilon}{8+3\varepsilon},~r=\tilde{p}=\frac{6p}{4p-3}.$

For any
$p\in(1,\frac{3}{2})$, $~r\in (3,6)$, we can choose admissible pairs $(q,r)\in \Lambda$ such that  $q\in (2,4)$,
provided  that  $\varepsilon >0$ sufficient small. Hence
$$I_2\lesssim
T^{-{\frac{\varepsilon}{4(4+\varepsilon)}}}\big(\|V\|_p\|\Psi_j\|_X\|U\|^2_X\big)\longrightarrow~0,~\text{~as~}
T\rightarrow \infty.
 $$

{\bf Step 2~}  Different with the proof of Step 1,  we perform an interior cut-off on space to get
\begin{align}&\Big\|\int_0^t\Delta^{-1}G'(t-s)\big((V*|U|^2)\eta_R\Psi_j\big)(s)ds\Big\|_{L^4(|t|\leq T,H^1_3(\R^3))}\nonumber\\
\lesssim& \|(V*|U|^2)\eta_R\Psi_j \|_{L^{4/3}(|t|\leq
T,H^1_{3/2})}\nonumber\\
\leq& \|(V*|U|^2)\eta_R \|_{L^2(|t|\leq T,L^{3})} \|\Psi_j
\|_{L^4(|t|\leq T,H^1_3)}\nonumber\\
&+ \|(V*|U|^2)\eta_R \|_{L^2(|t|\leq T,H^1_{2p})} \|\Psi_j \|_{L^4(|t|\leq T,L^{\tilde{p}})}\nonumber\\
\leq&\|\Psi_j \|_X\Big( \|(V*|U|^2)\eta_R \|_{L^2(|t|\leq
T,L^{3})}+\|(V*|U|^2)\eta_R \|_{L^2(|t|\leq
T,H^1_{2p})}\Big),\nonumber
\end{align} while
\begin{align}&\|(V*|U|^2)\eta_R \|_{L^2(|t|\leq
T,L^{3})}+\|(V*|U|^2)\eta_R \|_{L^2(|t|\leq T,H^1_{2p})}\nonumber\\
\leq&
\|V\|_p\|U\|^2_{L^4(|t|\leq T,L^{\tilde{p}})}+\|V\|_p\|U\|_{L^4(|t|\leq T,L^{\tilde{p}})}\|U\|_{L^4(|t|\leq T,H^1_3)}\nonumber\\
\leq& 2\|V \|_p \|U \|^2_{L^4(|t|\leq
T,H^1_{3})}<\infty.\nonumber\end{align}
Hence, we have
\begin{align} &~~~~\|(V*|U|^2)\eta_R \|_{L^2(|t|\leq T,L^{3})}+ \|(V*|U|^2)\eta_R
\|_{L^2(|t|\leq T,H^1_{2p})}\nonumber\\&=\|(V*|U|^2) \|_{L^2\big(
|t|\leq T ,L^{3}(|x|>R)\big)}+ \|(V*|U|^2)\eta_R
\|_{L^2\big(|t|\leq T ,H^1_{2p}(|x|>R)\big)}\nonumber\\
&\longrightarrow  ~0,~~\text{as~} R\longrightarrow~
+\infty.\end{align} \smallskip

 {\bf Step 3~} Observe  that
for each fixed $j\in\N$, $$\lim_{M\rightarrow\infty}
 \|{\cal T}_{U_0}(\xi_R\Psi_j) \|_{L^4(|t|\leq T,H^1_3(|x|>M))}=0.$$
In order to prove convergence uniformly on $j\in\N$, we take
advantage of the finite $\varepsilon$-cover property of compact set.

By Lemma \ref{lem3.1}, one has
$$ \|(V*|U|^2)\xi_R\Psi_j \|_{H^\delta(\R^3)}\leq \|(V*|U|^2)\Psi_j \|_{H^\delta(|x|\leq
2R)}\leq C.$$
This together with the Rellich-Kondrachov theorem implies that
$$\{(V*|U|^2)\Psi_j\}_{j=0}^\infty~\text{ is ~compact~ in ~}L^2(|x|\leq 2R).$$
This fact shows that  $\forall ~\varepsilon>0$, $\exists$ finite set
$\mathscr{A}=\{j_1,j_2,\cdots,j_{l_0}\}$ such that $~\forall~j\in\N$,
$\exists~ l\in\mathscr{A}$ satisfying
$$ \|(V*|U|^2)\xi_R\Psi_j-(V*|U|^2)\xi_R\Psi_l \|_{L^2}<\varepsilon.$$
Hence, as  $M$ is large enough,  we  obtain that by Lemma \ref{lem3.2}
\begin{align}&\Big\|\int_0^t\Delta^{-1}G'(t-s)\big((V*|U|^2)\xi_R\Psi_j\big)(s)ds\Big\|_{L^4(|t|\leq
T,H^1_3(|x|>M))}\nonumber\\
\leq&\Big\|\int_0^t\Delta^{-1}G'(t-s)\big((V*|U|^2)\xi_R\Psi_j-(V*|U|^2)\xi_R\Psi_l\big)ds\Big\|_{L^4(|t|\leq
T,H^1_3)}\nonumber\\
&+\Big\|\int_0^t\Delta^{-1}G'(t-s)\big((V*|U|^2)\xi_R\Psi_l\big)(s)ds\Big\|_{L^4(|t|\leq
T,H^1_3(|x|>M))}\nonumber\\
\leq&C(T)\sup_{t\in[-T,T]}\Big\|\int_0^t\Delta^{-1}G'(t-s)\big((V*|U|^2)\xi_R\Psi_j-(V*|U|^2)\xi_R\Psi_l\big)ds\Big\|_{H^{2-\epsilon}(\R^3)}+\varepsilon\nonumber\\
 \leq&
C(T) \sup_{t\in[-T,T]}\|(V*|U|^2)\xi_R\Psi_j-(V*|U|^2)\xi_R\Psi_l \|_{L^2(\R^3)}+\varepsilon\nonumber\\
\lesssim& ~\varepsilon\nonumber\end{align}

\smallskip
{ \bf Step 4~}\; After localizing $t$ and $x$ to bounded domain, we can use the following Arzela-Ascoli
compactness argument.

\begin{lemma}\label{lem3.3} A sequence $\{f_j\}_{j=0}^\infty$ in
$C([-T,T],H^1_3(|x|\leq M))$ has a convergent subsequence iff

 (i) for each $t\in [-T,T]$, the sequence $\{f_j(t)\}_{j=0}^\infty$ has
 a convergent subsequence in $H^1_3(|x|\leq M)$;

 (ii) the sequence $\{f_j\}_{j=0}^\infty$ is equicontinuous  on
 $[-T,T]$.
\end{lemma}

We now verify that
$$f_j(t)=-\int_0^t\Delta^{-1}G'(t-s)\big((V*|U|^2)\xi_R\Psi_j\big)(s)ds$$ satisfies
the two conditions of Lemma 3.3.

 By Lemma \ref{lem3.2}, we have
for all $t\in[-T,T]$
\begin{align}\label{3.2}&\Big\|\int_0^t\Delta^{-1}G'(t-s)\big((V*|U|^2)\xi_R\Psi_j\big)(s)ds\Big\|_{H^{2-\epsilon}(|x|\leq
M))}\nonumber\\\lesssim&  \|(V*|U|^2)\xi_R\Psi_j
\|_{L^2([-T,T],L^2)}\leq C.\end{align}
This, together with the Rellich-Kondrachov
theorem, implies that the sequence $\{f_j\}_{j=0}^\infty$ satisfies  (i) of in Lemma \ref{lem3.3}.

Next, we show that  the equicontinuity of the sequence
$\{f_j\}_{j=0}^\infty$ on
 $[-T,T]$.
\begin{align}& \|f_j(t+h)-f_j(t) \|_{H^1_3}\nonumber\\
=&\Big\|\int_0^{t+h}\Delta^{-1}G'(t+h-s)\big((V*|U|^2)\xi_R\Psi_j\big)(s)ds\nonumber\\
&-\int_0^t\Delta^{-1}G'(t-s)\big((V*|U|^2)\xi_R\Psi_j\big)(s)ds\Big\|_{H^1_3}\nonumber\\
\leq&\Big\|\int_0^t\big(\Delta^{-1}G'(t+h-s)-\Delta^{-1}G'(t-s)\big)\big((V*|U|^2)\xi_R\Psi_j\big)(s)ds\Big\|_{H^1_3}\nonumber\\
&+\Big\|\int_t^{t+h}\Delta^{-1}G'(t+h-s)\big((V*|U|^2)\xi_R\Psi_j\big)(s)ds\Big\|_{H^1_3}\nonumber\\
=&:J_1+J_2\nonumber\end{align}

 Let $U_j:=(V*|U|^2)\xi_R\Psi_j$.  By Lemma \ref{lem3.2} and  the
compactness as same as that in Step 3, it is derived that
\begin{align}\label{3.3}J_1=&\Big\|\int_0^t\Delta^{-1}G'(t-s)(G(h)-I)U_j(s)ds\Big\|_{H^1_3}\nonumber\\
\lesssim&\Big\|\int_0^t\Delta^{-1}G'(t-s)(G(h)-I)U_j(s)ds\Big\|_{H^{2-\epsilon}}~~~(\text{for some }\epsilon>0)\nonumber\\
\lesssim&  \|(G(h)-I)U_j \|_{L^2([-T,T],L^2)}\nonumber\\
\leq& \|(G(h)-I)(U_j-U_l) \|_{L^2([-T,T],L^2)}+ \|(G(h)-I)U_l\ \|_{L^2([-T,T],L^2)}\nonumber\\
\leq&2 \|U_j-U_l \|_{L^2([-T,T],L^2)}+ \|(G(h)-I)U_l
\|_{L^2([-T,T],L^2)}<\varepsilon\end{align}
uniformly on $j\in\N$ as
$|h|$ is small enough.  Combining the $L^p$-$L^{p'}$ estimate with the H\"{o}lder inequality, we deduce that
\begin{align}\label{3.4}J_2\leq&
\int_t^{t+h}|t+h-s|^{-\frac{1}{2}} \|V*|U|^2\Psi_j \|_{H^1_{3/2}}ds\nonumber\\
\leq&\Big(\int_t^{t+h}|t+h-s|^{-\frac{q'}{2}}ds\Big)^\frac{1}{q'}
\|V*|U|^2\Psi_j \|_{L^qH^1_{3/2}},\end{align} where
$\frac{1}{q}+\frac{1}{q'}=1.$

Let $1+\frac{2}{3}=\frac{1}{p}+\frac{3}{r}$. One easily verify
$3<r<\frac{9}{2}$ for any $1<p<\frac{3}{2}$. This allows us to
choose admissible pairs ~$(q,r)\in \Lambda$ such that
\begin{align} & \|V*|U|^2\Psi_j \|_{L^qH^1_{3/2}}\nonumber\\
\leq& \|V \|_p( \|U \|^2_{L^\infty L^r}
 \|\Psi_j \|_{L^qH^1_r}+ \|U \|_{L^\infty L^r} \|\Psi_j \|_{L^\infty
L^r} \|U \|_{L^qH^1_r})\nonumber\\
\leq&c \|V \|_p \|U \|^2_X \|\Psi_j \|_X\leq C\label{bod}\end{align}
and
\begin{align}\int_t^{t+h}|t+h-s|^{-\frac{q'}{2}}ds\rightarrow
0, ~~\text{as}~~h\rightarrow 0 .\label{decay}\end{align} Hence,
$J_2\rightarrow 0 $ uniformly on $j\in\N$ as $h\rightarrow 0$.

\medskip
{\bf Step 5~} A Cantor diagonalized process.

 For each $N\in\{1,2,3,\cdots\}$, we  choose a $T(N)$ in Step 1, then a
 $R(N)$ in Step 2 and then a $M(N)$ in Step 3 such that
\begin{align}\label{12}&\sup_{j\in\N} \|\mathcal{T}_{U_0}\Psi_j \|_{L^4(|t|>T(N),H^1_3(\R^3))}<\frac{1}{N},\\
&\sup_{j\in\N} \|\mathcal{T}_{U_0}(\eta_{R(N)}\Psi_j) \|_{L^4(|t|\leq T(N),H^1_3(\R^3))}<\frac{1}{N},\\
&\sup_{j\in\N} \|\mathcal{T}_{U_0}(\xi_{R(N)}\Psi_j)
\|_{L^4(|t|\leq T(N),H^1_3(|x|>M(N)))}<\frac{1}{N}.\end{align}
In this way, we  can choose inductively subsequence $\{\Psi_{j,N}\}$ of
$\{\Psi_{j,N-1}\},~N=1,2,\cdots $ with $\Psi_{j,0}=\Psi_j$, such
that $\{\mathcal{T}_{U_0}(\xi_{R(N)}\Psi_{j,N})\}$ converges in
$L^4(|t|\leq T(N),H^1_3(|x|\leq M(N)))$. Thus the subsequence
$\{\mathcal{T}_{U_0}(\xi_{R(N)}\Psi_{N,N})\}_{N=1}^\infty$
converges in $L^4(\R,H^1_3(\R^3))$. This completes the proof of
Lemma \ref{main}.

\section{Proof of Theorem \ref{th1}}

 In this  paper we still take advantage of the approach in
Kumlin \cite{Ku} to prove Theorem 1.1, and  by exploring sufficiently
compactness condition we give a more concise proof.
\begin{lemma}$^{\cite{HP}}$ \label{lem4.1} Let $H $ be Hilbert space,
$A_k:~H\mapsto H,~k=1,2,\cdots,$ be ananlytic
mappings, uniformly bounded on all compact set $D\subset H$. Also
assume that $A_ku\rightarrow Au$ as $k\rightarrow\infty$ for all
$u\in H$. Then the mapping $A:~H\rightarrow H$ is also analytic.
\end{lemma}

According to Theorem \ref{th2}, one has that $\mathcal
{N}(T):~U_0\mapsto U(T)$ is analytic from $H^1$ to $H^1$ for every
$T\in\R$. The wave operators $W_{\pm}$ and their inverses can be
represented as
\begin{align}\label{4.2}W_{\pm}=\lim_{T\rightarrow\pm\infty}\mathcal {N}(-T)G(T)\\
\label{4.3}W_{\pm}^{-1}=\lim_{T\rightarrow\pm\infty}G(-T)\mathcal
{N}(T)\end{align} Note that $\mathcal {N}(-T)G(T)$ and
$G(-T)\mathcal {N}(T)$ are analytic on $H^1$, and
$G(T)$ is an isometric on $H^1$, Lemma \ref{lem4.1}
implies that $W_{\pm},~W_{\pm}^{-1}$ and $S$ are analytic provided that
\begin{align}\label{4.4}\sup_{\Phi\in D}\sup_{T\in\R} \|\mathcal {N}(T)\Phi \|_{H^1}<\infty\end{align}
for all compact set $D\subset H^1$. In fact,
\begin{align} \|\mathcal {N}(T)\Phi \|_{H^1}&\leq  \|G(T)\Phi \|_{H^1}
+ \Big\|\int_0^T\Delta^{-1}G'(t-s)(V*|U(\Phi)|^2)U(\Phi)(s)ds \Big\|_{L^\infty
H^1}\nonumber\\&\leq \|\Phi \|_{H^1}+\|V\|_p \|U(\Phi)
\|^3_{L^4H^1_3}.
\end{align}Hence, it is enough to prove
\begin{align}\label{4.6}\sup_{\Phi\in D} \|U(\Phi) \|_{L^4H^1_3}<\infty.\end{align}
Now we give the proof of \eqref{4.6} briefly by using the finite $\varepsilon$-cover again and the definition of
Frech\'{e}t derivative.

 In fact,  As a direct result of the scattering theory, we have
$$ \|U(\Phi) \|_{L^4H^1_3}<\infty$$ for each $\Phi\in H^1$.
Hence, we  need to prove it is bounded uniformly on $\Phi\in D$.

 Since $D$ is a
compact subset of $H^1$, then for fixed $0<\varepsilon_0<1$, there exists a  finite set
$\mathscr{A}=\{\Phi_{l_1},\Phi_{l_2},\cdots,\Phi_{l_0}\}$ such that for
any $\Phi\in D$, there exists $\Phi_l\in\mathscr{A}$ satisfying
$$ \|\Phi-\Phi_l \|_{H^1}<\varepsilon_0.$$
Note that $U:~\Phi\rightarrow U(\Phi)$ is analytic from $H^1$ to
$L^4H^1_3$, we easily see that the  Frech\'{e}t derivative $U'(U_l)$ is a bounded operator
from $H^1$ to $L^4H^1_3$. This yields that
\begin{align} \|U(\Phi) \|_{L^4H^1_3}&\leq \|U(\Phi)-U(\Phi_l) \|_{L^4H^1_3}
+ \|U(\Phi)-U(\Phi_l) \|_{L^4H^1_3}\nonumber\\
&\leq
 \|U'(\Phi_l)(\Phi-\Phi_l) \|_{L^4H^1_3}+o(\varepsilon_0)+ \|U(\Phi_l) \|_{L^4H^1_3}\nonumber\\
&\leq
C_l \|\Phi-\Phi_l \|_{H^1}+o(\varepsilon_0)+ \|U(\Phi_l) \|_{L^4H^1_3}\nonumber\\
&\leq C_l\varepsilon_0+o(\varepsilon_0)+ \|U(\Phi_l)
\|_{L^4H^1_3}<C.\end{align}
This completes  the proof of Theorem
\ref{th1}.

\bigskip
\bigskip
 {\bf Acknowledgments:}\quad The authors thank the referees and the associated
editor for their invaluable comments and suggestions which helped
improve the paper greatly. The authors are grateful to Prof. W.
Strauss for his valuable suggestions.   C.Miao  was partially  supported
 by the NSF of China (No.10725102).

.
\begin{center}

\end{center}
\end{document}